\documentclass[12pt]{amsart}

\newtheorem{theorem}{Theorem}[section]
\newtheorem{lemma}[theorem]{Lemma}
\newtheorem{proposition}[theorem]{Proposition}

\newtheorem{definition}[theorem]{Definition}

\newtheorem{corollary}[theorem]{Corollary}

\newcommand{\Co}{\mbox{$\mathbb{C}$}}
\newcommand{\N}{\mbox{$\mathbb{N}$}}

\newcommand{\K}{\mbox{$\mathbb{K}$}}

\newcommand{\C}{\mbox{${\mathcal C}$}}
\newcommand{\D}{\mbox{${\mathcal D}$}}

\newcommand{\G}{\mbox{${\mathcal G}$}}

\newcommand{\Sy}{\mbox{${\mathcal S}$}}
\newcommand{\Li}{\mbox{${\mathcal L}$}}
\newcommand{\M}{\mbox{${\mathcal M}$}}
\newcommand{\be}{\mbox{${\mathbb B}$}}
\newcommand{\U}{\mbox{${\mathcal U}$}}

\begin{document}
\title[Multipliers and injectivity]{Multipliers
of operator spaces, and the injective envelope}

\vspace{30 mm}

\author{David P. Blecher* and Vern I. Paulsen*}
\address{Department of Mathematics\\University of Houston\\Houston,
TX 77204-3476 }
\email{dblecher@math.uh.edu  vern@math.uh.edu}
\date{August 27, 1999.  Revised January 21, 2000.}
\thanks{PORTIONS OF THIS WORK WERE PRESENTED AT THE
 CANADIAN OPERATOR THEORY AND OPERATOR ALGEBRAS
SYMPOSIUM, May 20 1999}
\thanks{* Supported by a grant from the NSF}

\maketitle

\vspace{40 mm}

\begin{abstract}
We study the injective envelope $I(X)$ of an operator space $X$,
showing amongst other things that it is a self-dual C$^*-$module.
We describe 
the diagonal corners of the injective envelope of the canonical 
operator system associated with $X$.
We prove that if $X$ is an operator $A-B$-bimodule, then $A$ and $B$ can be
represented completely contractively as subalgebras of these corners.
Thus, the operator algebras that can act on $X$ are determined by these
corners of $I(X)$ and consequently bimodules actions on $X$ extend naturally
to actions on $I(X)$.
These results give another characterization of the multiplier algebra of an 
operator space, which was introduced by the first author,
and a short proof of a recent characterization of operator modules,
and a related result. 
As another application, we extend Wittstock's module map extension
theorem, by showing that an operator $A-B$-bimodule is injective
as an operator $A-B$-bimodule if and only if it is injective as an
operator space. \end{abstract}

\pagebreak
\newpage

\section{Introduction.}  
In this paper we investigate some connections between the 
following topics: injectivity
of operator spaces, self-dual Hilbert C$^*-$modules (in the sense of
 Paschke);  completely contractive
actions of one operator space on another;
and the notion of a multiplier operator algebra
 of an operator space which was recently introduced by the first author.

The results and definitions follow a natural logical sequence, so 
we begin without further delay.   
We refer to \cite{Bsh} for additional information, and complementary
results, and to \cite{P,Pis} for background information on operator
spaces and completely bounded maps.

Recall that an operator space $X$ is
{\em injective}, if for any operator spaces $W \subset Z$, and any
completely bounded linear $T : W \rightarrow X$, there exists a
linear $\tilde{T} : Z \rightarrow X$ extending $T$, with 
$\Vert T \Vert_{cb} = \Vert \tilde{T} \Vert_{cb}$.  It has been 
known for a long time that $B(H)$ is an injective operator 
space (see \cite{P,Wi,P1,Arv}) for any Hilbert space $H$.
In 1983, Youngson  showed \cite{Yo} that
an injective operator space $X$ is a `corner' of a C$^*-$algebra $A$,
by which we mean that there exists projections $p,q$ in the multiplier 
algebra $\M(A)$ of $A$, such that $X = p A q$.   It is well known
\cite{Ri2} that this last condition
is equivalent to saying that such an $X$ is a Hilbert
C$^*-$module.  Then Hamana in 1985 (see notes in
\cite{Ham3}), 
and Ruan independently
a little later \cite{Rua}, 
showed that any operator space $X$ has an
operator space injective envelope $I(X)$.  To prove the 
existence of this envelope, one may follow
the classical construction for Banach space injective envelopes.
We will sketch the main idea: 
one begins by choosing any injective 
object $B$ containing $X$.  Then one considers the {\em $X$-projections}
on $B$, by which we mean completely contractive idempotent maps
on $B$ which fix $X$.  An idempotent map of course
is one such that $\Psi \circ \Psi = \Psi$.
There is a natural ordering on such maps,
and, with a little work one can show,
by a Zorn's lemma argument, that there is a minimal $X$-projection
$\Phi$.
The range of $\Phi$ in $B$ may be taken to be the injective 
envelope $I(X)$ of $X$, and one has $X \subset I(X) \subset B$.  Thus
one sees that $I(X)$ is the smallest injective space containing $X$.
As in the Banach space case (see \cite{La} for details and references),
one proves that $I(X)$ is an `essential'
and `rigid' extension of $X$.  The latter term, rigidity, means 
that the identity map is the only completely contractive
map on $I(X)$ extending the identity map on $X$.   

Note that if 
$\Sy$ is a linear subspace of $B(H)$ containing $I_H$ (for 
example, if $\Sy$ is a unital C$^*-$algebra), 
then one may choose $B =B(H)$ in the 
above.  Since $\Phi(I) = I$, it follows that $\Phi$ is 
completely positive \cite{P}.  A well-known theorem of 
Choi and Effros  \cite{CE} states that the range of a 
completely positive unital idempotent 
map on a C$^*-$algebra, is a 
C$^*-$algebra with respect to a certain multiplication.
Hence it follows that $I(\Sy)$ is a unital C$^*-$algebra\footnote{In
fact, $I(A)$ is a unital C$^*-$algebra even for a non-unital C$^*-$algebra 
$A$.  Since we have no good reference for this no doubt well
 known fact, we supply a proof later.}.     

Hamana gives a different construction of
$I(X)$ to the one outlined above, which allows one to 
prove something a little stronger.  Since this construction will
be important for us, we will outline some of the ideas.
By the  method popularized by the second author
(see \cite{P} Lemma 7.1), we may 
embed the operator space $X$ in a
canonical unital operator system\footnote{An operator system is
a selfadjoint linear subspace of $B(H)$ containing $I_H$.}
$$ \Sy(X)  =\left(
\begin{array}{cc}
\Co & X \\
X^* & \Co
\end{array}
\right)  \; \; \; . 
$$

If one forms the injective envelope
$I(\Sy(X))$, it will be a unital C$^*-$algebra, by the 
aforementioned argument using the Choi-Effros result.
Indeed, since the minimal $\Sy(X)$-projection fixes the 
C$^*-$algebra $\Co \oplus \Co$ which is the diagonal of
$\Sy(X)$, it follows immediately (for example by Lemma
\ref{chtr} below, although this is not necessary), 
that the following elements of $\Sy(X)$ are two 
selfadjoint projections with sum 1 in  the C$^*-$algebra 
$I(\Sy(X))$:
$$ p = \left(
\begin{array}{cc}
1 & 0 \\
0 & 0 
\end{array}
\right)
\; \; \; \; , \; \; \; \; q \; = \; \left(
\begin{array}{cc}
0 &  0 \\
0 & 1 
\end{array}
\right)
\; \; . $$      
Therefore, with respect to $p$ and $q$, we may decompose $I(\Sy(X))$ to 
write it as consisting of $2 \times 2$ matrices.  Hamana shows
 that $p I(\Sy(X)) q,$
the 1-2 corner of $I(\Sy(X))$, is the injective envelope 
of $X$.  This recovers and strengthens Youngson's
result.
The four corners of $I(\Sy(X))$ we will name:
$$I(\Sy(X)) = \left(
\begin{array}{cc}
I_{11} & I(X) \\
I(X)^* & I_{22} 
\end{array}
\right) \; \; .
$$
It is clear that $I_{11}$ and 
$I_{22}$ are also  injective 
C$^*-$algebras.

We will write $j$ for the canonical inclusion of 
$X$ inside $I(X)$.   
 
Writing $I = I(X)$ for a moment, we define a subset of 
$I(\Sy(X))$ by 
$$\Li(I(X)) = \left(
\begin{array}{cc}
I I^* & I \\
I^* & I^*I
\end{array}
\right) \; \; .
$$
where $I I^*$ for example, is the closed span in $I_{11}$
of terms $x y^*$ with $x,y \in I$.  
Henceforth, we will reserve the letters $\C(X)$ and $\D(X)$
for $I I^*$ and $I^*I$ respectively.
 When $X$ is understood
we will simply write $\C$ and $\D$.  Thus $\C \subset
I_{11} \; , \D \subset I_{22}$.  Notice also that 
$\Li(I(X))$ coincides with the smallest closed 2-sided ideal in 
$I(\Sy(X))$ containing the copy of $I(X)$ in the 1-2-corner 
of $I(\Sy(X))$,    
and this fact will be used below.

Clearly 
$\Li(I(X))$ is a C$^*-$subalgebra of $I(\Sy(X))$.  From this,
or otherwise, it is easy to see that $I(X)$ is a $\C-\D-$bimodule
which, with respect to the natural $\C$- and $\D$-valued inner
products given by $xy^*$ and $x^*y$, is a `strong Morita 
equivalence $\C-\D-$bimodule' in the language of Rieffel.
Sometimes this is also referred to as a `$\C-\D-$imprimitivity
bimodule'.  It follows by basic C$^*-$algebraic Morita theory
(see \cite{Ri2,L2} say), that 
$\C \cong \K(I(X))$ as C$^*-$algebras, where $\K(I(X))$ is the 
so called `imprimitivity C$^*-$algebra' of the right $\D$-module
$X$.  Also, $\M(\C) \cong \be_{\D}(I(X))$, the adjointable
$\D$-module maps on $I(X)$.  In fact we shall see that
the left multiplier algebra $LM(\C)$ (which by a result of Lin
\cite{Lin} may be identified with the space
$B_{\D}(I(X))$ of  bounded right
$\D$-module maps
on $I(X)$) coincides with $\M(\C)$.   Also $\Li(I(X))$ may be 
identified with the `linking C$^*-$algebra' \cite{Ri2}
of the bimodule $I(X)$.  Although this will not be explicitly 
used below, it is a useful perspective.

We shall write $\Sy_0(X)$ for the subspace of $\Sy(X)$
consisting of those
elements with $0$'s on the main diagonal.
 
\begin{proposition}
\label{ess}
For any operator space $X$, we have that 
$J = I(\Sy(X)) \Sy_0(X) I(\Sy(X))$ is an essential ideal in 
$I(\Sy(X))$.
\end{proposition}

\begin{proof}  Suppose that $K$ is an ideal in $I(\Sy(X))$
whose intersection with $J$ is zero.  Let $\pi :
I(\Sy(X)) \rightarrow I(\Sy(X))/K$.  Since $\pi$ is 1-1 on
$J$, it is completely isometric on $\Sy_0(X)$.  Let $\Phi$ be the 
restriction of $\pi$ to $\Sy(X)$.  By Lemma 7.1 in 
\cite{P}, $\Phi$ is a complete order injection.    

Extend the map from $\pi(\Sy(X)) \rightarrow \Sy(X)$ which
is  the inverse of $\Phi$, to a map $\gamma :
I(\Sy(X))/K \rightarrow I(\Sy(X))$.  Since 
$\gamma \circ \Phi = Id_{\Sy(X)}$, it follows by rigidity
that $\gamma \circ \pi = Id_{I(\Sy(X))}$.  Thus 
$K = (0)$.
\end{proof}

\begin{corollary}
\label{es2}  For any operator space $X$, 
we have that $\Li(I(X))$ is an essential ideal
in $I(\Sy(X))$, and that $\C = \C(X)$ is an essential ideal
in $I_{11}$.  Thus we have the following canonical inclusions 
of C$^*-$algebras
$$\C \subset I_{11} \subset M(\C) \; .$$
If $\C$ is represented faithfully and nondegenerately
on a Hilbert space $H$, then this string may be regarded as 
inclusions of subalgebras of $B(H)$.  Similar assertions
hold for $\D(X) \subset I_{22}$.
\end{corollary}

\begin{proof}  Clearly $J \subset \Li(I(X))
\subset I(\Sy(X))$, where $J$ is as in Proposition
\ref{ess}.  Thus $\Li(I(X))$ is an essential ideal
in $I(\Sy(X))$, by that Proposition.
To see the second assertion, notice
that if $t \in I_{11}$ and $t \C = 0$, then 
$t \C t^* = 0$, which implies that $t z z^* t^* = t z = 0$ 
for all $z \in I(X)$.  
It follows immediately that 
$(t \oplus 0) \Li(I(X)) = 0$.  Hence by the first assertion,
$t = 0$.  Thus $\C$ is an essential ideal
in $I_{11}$.  

That $I_{11} \subset M(\C)$ follows from 
the universal property of the multiplier algebra, namely 
that $M(\C)$ contains a copy of any C$^*-$algebra containing
$\C$ as an essential ideal (\cite{L2} Chapter 2, say).   
\end{proof}

\vspace{5 mm}

Since $\C$ is essential in $I_{11}$, it follows from 
\cite{FP} Theorem 2.12,
 that $M(\C) \subset I_{11}$.  Thus  
in fact $I_{11} = M(\C) \cong \be_{\D}(I(X))$.   However we will
 deduce all this in a self-contained way,
from some machinery we develop next:  

\begin{corollary}
\label{es3}  If $t \in I_{11}$, and if $t x = 0$ for all 
$x \in X$, then $t = 0$.
\end{corollary}

\begin{proof}  Without loss of generality, we may assume that
 $\Vert t \Vert \leq 1$.  By replacing
$t$ by $t^* t$, we may also suppose that $0 \leq t \leq 1$,
where this last `1' is the identity of
$I_{11}$.  Let $p = 1-t$.  Define $\phi(z) = p z$, for $z \in I(X)$.
Since $\phi(x) = x$ for $x \in X$, we obtain $\phi = Id$,
by rigidity.   Thus $t I(X) = 0$, which (we showed 
in the proof of Corollary \ref{es2}) implies that
$t = 0$.
\end{proof}
  
\begin{definition}  Let $X$ be an operator space.  We define
the {\em left multiplier operator algebra} of $X$ to be
$IM_l(X) = \{ T \in I_{11} : TX \subset X \}$.  We define the
left {\em multiplier C$^*-$algebra} of $X$ to be
$IM^*_l(X) = \{ T \in IM_l(X) : T^* \in IM_l(X) \}$.
\end{definition}

Note that $IM^*_l(X)$ is a C$^*-$algebra, whereas
$IM_l(X)$ is a unital nonselfadjoint operator algebra
in general.  
There is a similar definition for right multipliers.  

We shall soon see that these
multiplier algebras coincide with the ones introduced in
\cite{Bsh} \S 4.  These simultaneously generalize the 
common operator algebras associated with Hilbert C$^*-$modules,
the multiplier function algebras of a Banach space introduced
by Alfsen and Effros \cite{AE} (who only considered the real
scalar case - see \cite{Be} for the complex case), and the 
multiplier algebras of a nonselfadjoint operator algebra with
c.a.i..

\begin{definition}  Let $H, K$ be Hilbert spaces.
A {\em multiplication situation} on $H \oplus K$
consists of 
three concrete operator spaces $X,Y,Z$ such that
$Y \subset B(H), \; X \subset B(K,H),$
and $Z \subset B(K)$, and such that 
$YX \subset X$ and $XZ \subset X $.  
\end{definition}

We will need the following elementary lemma (c.f. \cite{P}
Ex. 4.2-4.5):

\begin{lemma}
\label{chtr}  (Choi \cite{ch} )
 Suppose that $\phi : A \rightarrow B$ is a
completely positive map between C$^*-$algebras, with $\phi(1) = 1$.
Suppose that there is a C$^*-$subalgebra $N$ of $A$ with 
$1_A \in N$, such that $\pi = \phi_{|_N}$ is a *-homomorphism.
Then $\phi$ is a `N-bimodule map'.  That is,
$$ \phi(an) = \phi(a) \pi(n) \; \; \; \; \text{and} \; \; \; \; 
\phi(na) = \pi(n) \phi(a) \; \; $$
for all $a \in A, n \in N$.
\end{lemma}

We will refer to the following as the `multiplication 
theorem':

\begin{theorem}
\label{op}  
\begin{itemize}
\item [(i)]  Suppose that $X$ is an operator space, and that 
$I_{11}, I_{22}$ are the diagonal corners of $I(\Sy(X))$, as usual.
If $Y,Z$ are two operator spaces such that 
$X, Y,Z$ form a multiplication situation on $H \oplus K$, as above, 
then 
there exists unique completely contractive linear maps
$\theta : Y \rightarrow IM_l(X) \; $ and $
 \pi : Z \rightarrow IM_r(X)$ 
such that $\theta(y) j(x) = j(yx)$, and 
$j(x) \pi(z) = j(xz)$, for all 
$x \in X, y \in Y , z \in Z$.
\item [(ii)]  If in addition,
$Y$ is a subalgebra (resp. *-subalgebra) of $B(H)$, 
then $\theta$ is also a homomomorphism (resp. *-homomomorphism into 
$IM_l^*(X)$).   Similarly for $Z$.   
\end{itemize} 
\end{theorem}

\begin{proof}  By \cite{P} Lemma 7.1, we have 
a completely order isomorphic copy $\G$
of $\Sy(X)$ inside
$B(H \oplus K)$.  By \cite{Ham} Corollary 4.2, there 
exists a surjective *-homomorphism from the C$^*$-subalgebra
$C^*(\G)$ of $B(H \oplus K)$ generated by $\G$, onto the
C$^*-$envelope $C^*_e(\Sy(X))$, 
which fixes the copies of $\Sy(X)$.  
Let $\phi : B(H \oplus K) \rightarrow I(\Sy(X))$ be a
completely positive map extending the *-homomorphism.
Since $\phi$ fixes the diagonal scalars
$\Co \oplus \Co$, it follows by
\ref{chtr} (this is a common argument), 
that $\phi$ decomposes as a $2 \times 2$ matrix of maps,
each corner map defined on the corresponding
 `corner' of $B(H \oplus K)$.
In particular, we have
$$
\phi(\left[ \begin{array}{ccl}
y & x \\
0   & z 
\end{array}
\right]  ) \; =
\left[ \begin{array}{ccl}
\theta(y) & j(x) \\
0   & \pi(z) \end{array}
\right] \; \; ,
$$
for a map $\theta : Y \rightarrow I_{11}$, 
$\pi : Z \rightarrow I_{22}$, and for all 
$x \in X, y \in Y, z \in Z$.  
By the previous lemma, $\phi$ is a `$C^*(\G)$-bimodule map'.
Hence for $y \in Y , x \in X$ we have:  
$$
\left[ \begin{array}{ccl}
0 & j(yx) \\
0   & 0 
\end{array}
\right]   = \phi \left(
\left[ \begin{array}{ccl}
y & 0 \\
0   & 0  
\end{array}
\right] 
\left[ \begin{array}{ccl}
0 & x \\
0   & 0 
\end{array}
\right]  \right) =
\phi \left(
\left[ \begin{array}{ccl}
y & 0 \\
0   & 0  
\end{array}
\right] \right) \left[ \begin{array}{ccl}
0 & jx \\
0   & 0 
\end{array}
\right]  
 = \left[ \begin{array}{ccl}
\theta(y) & 0 \\
0   & 0
\end{array}
\right]  \; \left[ \begin{array}{ccl}
0 & jx \\
0   & 0
\end{array}
\right]  $$
Thus $j(y x) = \theta(y) j(x)$.  
The uniqueness of $\theta$ follows from Corollary \ref{es3}.
Similarly for $\pi$. 
 
For (ii), note that for $y_1,y_2 \in Y$ and $x \in X$
we have $(y_1 y_2) x = 
\theta(y_1 y_2) x = \theta(y_1) \theta(y_2) x$.
Now use Corollary \ref{es3} to conclude that 
$\theta$ is a homomorphism.
The last assertion follows from the fact that a contractive 
representation of a C$^*-$algebra is a *-homomorphism.
 \end{proof}

\vspace{4 mm}

We recall that a right C$^*-$module $Z$ over $D$ is `self-dual', if it 
satisfies the equivalent of the Riesz representation theorem for 
Hilbert spaces, namely that every $f \in B_D(Z,D)$ is given by the 
inner product with a fixed $z \in Z$.

\begin{corollary}
\label{co3}
If $X$ is an operator space, then
\begin{itemize}
\item [(i)]  $I_{11} = M(\C(X)) = LM(\C(X)) = RM(\C(X)) = 
QM(\C(X)) = I(\C(X))$.
Thus $I_{11} \cong \be_{\D}(I(X))
 = B_{\D}(I(X))$.  Similarly,
$I_{22} = M(\D) = LM(\D) = I(\D)$.     
\item [(ii)]  $I(X)$ is a self-dual right C$^*-$module 
(over $\D$ or over $\M(\D)$).  Similarly, it is a self-dual 
left C$^*-$module.
\item [(iii)]  $I(\Sy(X))$ is the multiplier C$^*-$algebra
of $\Li(I(X))$.  Also, $I(\Sy(X))$ is the injective 
envelope of $\Li(I(X))$.
\end{itemize} \end{corollary} 

\begin{proof}  (i): 
Represent $\Li(I(X))$ nondegenerately on a Hilbert space $H \oplus K$.  
We obtain a multiplication situation on $H \oplus K$ given by the 
actions of $LM(\C)$ and $RM(\D)$ on $I(X)$.    
By the multiplication theorem, we get a completely contractive 
homomomorphism $\theta : LM(\C) \rightarrow I_{11}$ 
such that $\theta(T) x = T x$, for all $x \in I(X)$ and $T \in LM(\C)$.   
Hence $T = \theta(T) \in I_{11}$.  Thus 
by Corollary \ref{es2}, we have
$I_{11} = M(\C(X)) = LM(\C(X))$.  By Lin's theorem \cite{Lin},
the latter C$^*-$algebra may be identified with $B_{\D}(I(X))$.
Taking adjoints gives
$RM(\C(X)) = M(\C(X))$.  As noted in \cite{Lin2},
 (ii) below implies that 
$M(\C(X)) = QM(\C(X))$.   Finally, note that we have the following
C$^*-$subalgebras: 
  $\C \subset I_{11} = M(\C) \subset I(\C)$.
The last inclusion is a fact from \cite{FP} which is 
reproved at the end of our paper.  By the injectivity
of $I_{11}$ there exists a completely contractive projection
$\Phi$
from $I(\C)$ onto $I_{11}$; since $\Phi$ fixes $\C$ we see
by rigidity that $\Phi = Id$ and $I_{11} = I(\C)$.

(ii):  We apply (i) but with $X$ replaced by the 
right C$^*-$module
sum $I(X) \oplus_c M(\D)$.  Clearly the latter is a full
C$^*-$module over $M(\D)$, and it is injective since it is
a corner of $I(\Sy(X))$.     
Using (i) and Corollary 7.8 in \cite{Bsh} if necessary,
we obtain that
$\be_{M(\D)}(I(X) \oplus_c M(\D)) = B_{M(\D)}(I(X) \oplus_c M(\D))$,
from which it is easy to see that $I(X)$ is a 
self-dual $M(\D)$-module.  However $B_{\D}(I(X),\D)
= B_{M(\D)}(I(X),M(\D))$ by Cohen's factorization theorem.

(iii):  It is well known that for any full right 
C$^*-$module $Z$ over a C$^*-$algebra $B$ say, we have that 
$\K_{\D}(Z \oplus_c \D)$ is a copy of
the linking C$^*-$algebra of $Z$.
Also $\be_{\D}(Z \oplus_c \D) = M(\K_{\D}(Z \oplus_c \D))$.
Thus we have $M(\Li(I(X))) = \be_{\D}(I(X) \oplus \D)$.
This of course splits into four corners. 
The 1-1 corner is $\be_{\D}(I(X)) = I_{11}$ by (i).
The 2-1 corner is $\be_{\D}(I(X),\D) \cong I(X)$ by (ii).
The 2-2 corner is $\be_{\D}(\D,\D) = M(\D) = I_{22}$
by (i).  This gives the first result.  The second 
follows just as $I(\C) = I_{11}$ in (i).
\end{proof}

\begin{theorem}
\label{Ieq}  If $X$ is an operator space, then
\begin{itemize}
\item [(i)]  $IM_l(X)$ and $IM^*_l(X)$ are completely isometrically
isomorphic to the left multiplier algebras 
$M_l(X)$ and $\be_l(X)$  defined in 
\cite{Bsh}.
\item [(ii)]  $IM^*_l(X)$ is isometrically isomorphic to 
a closed subalgebra of $B(X)$.
Also, $IM^*_l(X)$ is completely isometrically isomorphic to 
a closed subalgebra of $B_l(X)$ or $CB_l(X)$.
\end{itemize} \end{theorem}

Before we prove this, we define $B_l(X)$ and $CB_l(X)$ for
an operator space $X$.  Namely $B_l(X) = B(X)$ but with matrix
norms 
$$\Vert [T_{ij}] \Vert^l_n = \sup \{ \Vert 
[\sum_{k=1}^n T_{ik}(x_k) ] \Vert_{C_n(X)} \; :
\; x \in BALL(C_n(X)) \} \; . $$
(Here $C_n(X) = X^n$, but with the operator space structure one gets
by identifying $C_n(X)$ with the `first column' of the operator
space $M_n(X)$.) 
With these norms $B_l(X)$ is not a matrix normed space in the
traditional sense.  However $M_n(B_l(X))$ is a unital Banach 
algebra.  Similarly one defines $CB_l(X)$, the only difference
being that one replaces $x_k$ in the expression above by
$x_{(k,p),q}$.  Here $[ x_{(k,p),q} ]$ is a matrix indexed on
rows by $(k,p)$ and on columns by $q$.  Again 
$M_n(CB_l(X))$ is a unital Banach 
algebra. 

\begin{proof}  (i):  By \ref{co3}, any $T \in IM_l(X)$ 
may be viewed as a bounded, and hence
completely bounded, module map on $I(X)$.
We obtain a canonical sequence of 
completely contractive homomorphisms
$$ IM_l(X) \overset{\rho}{\rightarrow}  
M_l(X) \rightarrow CB(X) \; $$
given by restriction of domain.  By Corollary \ref{es3},
these homomorphisms are 1-1.

On the other hand, we have a multiplication situation
 given by the 
action of $M_l(X)$ on $X$.  Hence, by the multiplication 
theorem, there exists a 
completely contractive homomorphism $\theta : 
M_l(X) \rightarrow I_{11}$ such that 
$\theta(T) x = Tx$ for all $T \in M_l(X), x \in X$.
Thus $\rho \theta = Id$.  Thus $\rho$ is onto, and
since $\rho$ is 1-1 we obtain that
$IM_l(X) \cong M_l(X)$ completely isometrically
and as operator algebras.
  
(ii) The first statement follows
from a result which may be found in
\cite{To} Proposition 1.1 or
\cite{B1} Corollary 1, 
which asserts that any
contractive homomorphism from a C$^*$-algebra into a
Banach algebra, is a $*-$homomorphism onto
its range, which is a C$^*$-algebra (with the norm and algebra
structure inherited from the Banach algebra).
Hence the canonical homomorphism 
$\be_l(X) \rightarrow CB(X)$ (or $\be_l(X) \rightarrow B(X)$)
is a $*-$homomorphism onto a C$^*$-algebra.  
Since the homomorphism is 1-1 
it is therefore isometric.  

The second statement
follows by considering 
the following canonical isometric inclusions
$$M_n(\be_l(X)) \subset M_n(\be_{\D}(I(X))) \subset 
M_n(B_l(I(X))) \; , $$
the last inclusion following from section 3 of \cite{Mg}.
Thus by restriction of domain, we get a contractive unital 
homomorphism $M_n(\be_l(X)) \rightarrow M_n(B_l(X))$.  
Now we can apply \cite{B1} Corollary 1 to deduce this
last homomorphism is an isometry.

A similar argument works for $CB_l$.   
\end{proof}

\vspace{4 mm}

It follows from (i) and a result in 
\cite{Bsh}, that 
the subalgebra of $B(X)$ or $CB(X)$ corresponding to
$IM_l^*(X)$ by (ii) above, is the 
C$^*-$algebra of (left) {\em adjointable operators}
 ${\mathcal A}_l(X)$ on $X$.

We do not know whether $IM_l^*(X)$ is completely 
isometrically contained inside $CB(X)$ in general.  
However it is not hard to find examples 
showing that the canonical contraction
$M_l(X) \rightarrow B(X)$  (or into $CB(X)$)
is not an isometry in general (see \cite{Bsh}).
This shows that if $T \in \be_{\D}(I(X))$, 
with $T(X) \subset X$,
then one cannot expect anything like
 $\Vert T \Vert_{cb} = \Vert T_{|_X} \Vert_{cb}$.

Finally we remark
 that if $IM_l(X) = \Co$, then it follows from the multiplication
theorem that for
any linear complete isometry $i : X \rightarrow B(K,H)$ such 
that $[i(X)K]^{\bar{}} = H$, we have that scalar multiples
of $I_H$ are  the only operators
 $T \in B(H)$ such that $T i(X) \subset i(X)$.

\vspace{5 mm}

\section{Applications.}

We recall from \cite{Bsh} that an {\em oplication} of an 
operator space $Y$ on an operator space $X$, is a 
bilinear map $\circ : Y \times X \rightarrow X$, such that 
\begin{itemize}
\item [(1)]  $\Vert y \circ x \Vert_n
\leq \Vert y \Vert_n \Vert x \Vert_n$, for all
$n \in \N , x \in M_n(X) , y \in M_n(Y)$,
\item [(2)]  There is an element $e \in Y_1$ such that
$e \circ x = x$ for all $x \in X$.
\end{itemize}

In (1), $y \circ x$ is computed by the usual rule for 
multiplying matrices.

\vspace{5 mm}
 
 We will apply the multiplication theorem from
 \S 1 to give a proof of the `oplication theorem' from 
\cite{Bsh} \S 5:

\begin{theorem}
\label{gen}  Suppose that $Y , X$ are operator spaces,
 and suppose that
$\circ : Y \times X \rightarrow X$ is an oplication, with
`identity' $e \in Y$.
Then there exists a unique completely contractive linear map
$\theta : Y \rightarrow M_l(X)$ such that
$y \circ x = \theta(y) x$ , for all $y \in Y, x \in X$.
Also $\theta(e) = 1$.  Moreover,  if $Y$ is, in addition,
an algebra with identity $e$, then $\theta$ is
a homomorphism if and only if $\circ$ is a module
action.  On the other hand,
if $Y$ is a C$^*-$algebra (or operator
system) with identity $e$,
then $\theta$ has range inside $\be_l(X)$, and is
completely positive and *-linear.
\end{theorem}

\begin{proof}  As in \cite{Bsh} the difficult part is to prove 
the first statement, the existence of 
$\theta$.  (The uniqueness follows from \ref{es3}.)
 Indeed it suffices
to find Hilbert spaces $H,K$, a complete isometry $\Phi : 
X \rightarrow B(K,H)$, and a linear complete contraction
$\theta : Y \rightarrow B(H)$, such that $\theta(e) = I$, 
and such that $\Phi(y  \circ x)
= \theta(y) \Phi(x)$, for all $y \in Y, x \in X$.  For in
that case $\theta(Y), \Phi(X) , \Co$ would form a 
multiplication situation on $H \oplus K$.  Now use 
the multiplication theorem above, together with the fact that
$IM_l(X) \cong M_l(X)$, to obtain the existence of $\theta$.

The existence of such $\Phi$ etc., follows easily from Le
Merdy's proof of the `BRS' theorem (\cite{LM} 3.3).  Namely,
first suppose that $X \subset B(K)$.   By the `multilinear 
Stinespring' theorem of \cite{CS,PS}, we may write
$y \circ x = \beta_1(y) \alpha_1(x)$, for completely 
contractive maps $\alpha_1 : X \rightarrow B(K,H_1)$ and 
$\beta_1 : Y \rightarrow B(H_1,K)$.  Similarly,
$\alpha_1(y \circ x) = \beta_2(y) \alpha_2(x)$, 
where now $\alpha_2 : X \rightarrow B(K,H_2)$ say.
Inductively we obtain, $\alpha_{k+1}(y \circ x) = \beta_k(y)
\alpha_k(x)$, where $\alpha_k : X \rightarrow B(K,H_k)$, say.
Since $e \circ x = x$, we see that each $\alpha_k$ is a 
complete isometry.
Let $W = X \otimes K$, and define $f_k : W \rightarrow H_k$
by $f_k(x \otimes \zeta) = \alpha_k(x)(\zeta)$.  It is easy to 
check that $\Vert f_k(w) \Vert_{H_k} \leq \Vert f_{k+1}(w) 
\Vert_{H_{k+1}}$, for $w \in W$, as in \cite{LM}.  Hence 
(by the parallelogram law if necessary) it is clear that 
$\lim_k \Vert f_k(\cdot) \Vert_{H_k}$ defines a seminorm 
which gives rise to a Hilbert space norm (on the quotient
of $W$ by the nullspace of the seminorm).   Write this 
resulting Hilbert space as $H$.  There is an obvious map
$\theta : Y \rightarrow L(H)$, given by $\theta(y)
([x \otimes 
\zeta]) = [(y \circ x) \otimes \zeta]$.  It is easy 
to see that this is completely contractive as in 
\cite{LM}.  The map $ \Phi : X \rightarrow B(K,H)$ given 
by $\Phi(x)(\zeta) = [x \otimes \zeta]$ is clearly a 
complete contraction, too.   
On the other hand, for any $\zeta \in K_1$, we have 
$$\Vert \Phi(x) \Vert \;  \geq \; \Vert \Phi(x)(\zeta) \Vert
\; = \; \lim_k \Vert \alpha_k(x)(\zeta) \Vert \; \geq \; 
\Vert \alpha_1(x)(\zeta) \Vert \; \; . $$
Thus $\Vert \Phi(x) \Vert \geq \Vert \alpha_1(x) \Vert =
\Vert x \Vert$, showing
that $\Phi$ is an isometry.  A similar argument shows that 
$\Phi$ is a complete isometry.   
\end{proof}

\vspace{5 mm}

As shown in \cite{Bsh}, this theorem has very many consequences,
containing as special cases, the `BRS' theorem \cite{BRS}, 
and many other results.   
The original proof of the above theorem in \cite{Bsh} was 
much more difficult.  The first author plans to replace that
proof with another, which quotes part of the
proof above, but uses
basic facts about `left order bounded operators' instead of the 
multiplication theorem.  The version in \cite{Bsh} also does
not give the fact that $IM_l(X) = M_l(X)$.

A left {\em operator module}
$X$ over a unital operator algebra $A$,
 is an operator space which is also a unitary left $A$-module
(unitary means that $1 x = x$ for all $x \in X$), such that 
$\Vert a x \Vert \leq \Vert a \Vert \Vert x \Vert$ for all
matrices $a$ with entries in $A$ and $x$ with entries in $X$.
A similar definition holds for right modules or bimodules.

One may deduce from the oplication 
theorem the following refinement of the
Christensen-Effros-Sinclair representation theorem
for operator modules \cite{CES}.  We will in fact give an 
independent proof (which is the essentially the same as the
proof above, except that we use the original
Christensen-Effros-Sinclair theorem (which is quite simple)
instead of the Le Merdy argument).   In any case note that
the method 
shows that the $A-B-$action on an operator $A-B$-bimodule
$X$, may be extended to an action on $I(X)$, making 
$I(X)$ an operator $A-B$-bimodule.    

\begin{theorem}
\label{nonchar}
Suppose that $A$ and $B$ are unital 
operator algebras, and that 
$X$ is an operator space and a unitary
$A-B$-bimodule with respect to a bimodule action 
$m : A \times X \times B \rightarrow X$.  
The following
are equivalent:
\begin{itemize}
\item [(i)]  $X$ is an operator $A-B$-bimodule.
\item [(ii)]  There exist Hilbert spaces $H$ and $K$,
and a linear complete isometry $\Phi : X \rightarrow B(K,H)$
and completely contractive unital homomorphisms 
$\theta : A \rightarrow B(H)$ and $\pi : B \rightarrow 
B(K)$, such that
$\Phi(m(a, x,b)) = \theta(a) \Phi(x) \pi(b)$, for all $a \in A,
x \in X$ and $b \in B$. 
\item [(iii)]  There exists unique completely contractive unital
homomorphisms
$\theta : A \rightarrow M_l(X)$ and 
$\pi : B \rightarrow M_r(X)$
such that $\theta(a) x \pi(b) = m(a, x,b )$
for all $a \in A, x \in X$ and $b \in B$.
\end{itemize}
\end{theorem}

\begin{proof}  That (i) is equivalent to (ii) is
a restatement of the Christensen-Effros-Sinclair 
representation theorem.  

$(iii) \implies (ii)$: Obvious.

$(ii) \implies (iii)$: $\theta(A), \Phi(X) , 
\pi(B)$ form a multiplication situation on $H \oplus K$.
The result then 
follows by the multiplication theorem. 
\end{proof}

\vspace{4 mm}

We will refer to a triple $(\Phi, \theta,\pi)$ as 
in (ii) above as a 
{\em CES representation} of the bimodule $X$.  
We will call it a {\em faithful CES representation} if
$\theta$ and $\pi$ are also completely isometric.
It is always possible, by an obvious direct sum trick,
to find a faithful CES representation for an 
operator bimodule.
For any CES representation of $X$ we obtain an
`upper triangular $2 \times 2$ operator
algebra', namely
$$ \U(X) = \left[ \begin{array}{ccl}
\theta(A) & \Phi(X) \\
0   & \pi(B)
\end{array}
\right] \; \; ,
 $$
We will write $\U_e(X)$ for this algebra in the case that 
$(j,\theta_e,\pi_e)$ is the representation
in (iii) above, into the multiplier 
algebras.   Thus $\U_e(X) \subset I(\Sy(X))$.

If $A$ and $B$ are C$^*-$algebras, and if we take 
faithful
CES representations and form $\U(X)$, then 
it is clear from \cite{Su} that $\U(X)$ as an abstract
operator algebra (i.e. as an algebra and an operator 
space), is independent of the particular faithful
CES representation.  This is probably not true 
if $A, B$ are nonselfadjoint.
However in either case, we can easily
 see that the triple $(j,\theta_e,\pi_e)$
given in (iii) above, is the 
`smallest' CES representation of $X$:  

\begin{corollary}
\label{min}  Suppose that $A$ and $B$ are unital operator 
algebras, and that $X$ is an $A-B$-operator bimodule. 
Let $(\Phi, \theta,\pi)$ be a CES representation
of $X$, and let $\U(X)$ be the corresponding 
upper triangular $2 \times 2$ operator
algebra.  Then there is a canonical 
completely contractive 
unital homomorphism $\phi : \U(X) \rightarrow \U_e(X)$,
which takes each corner of $\U(X)$ into the same corner
for $\U_e(X)$.   Indeed $\phi$ induces completely contractive
unital homomorphisms $\rho$ and $\sigma$, from 
$\theta(A)$ and $\pi(B)$, to $\theta_e(A)$ and 
$\pi_e(B)$ respectively, such that $\rho \circ \theta
= \theta_e$ and $\sigma \circ \pi = \pi_e$.   
\end{corollary}

\begin{proof}  This follows immediately from the proof
of Theorem \ref{op}, and of the implication `$(ii) 
\implies (iii)$' of \ref{nonchar}.  The $\phi$ is as in
the proof of Theorem \ref{op}, which is easily seen 
to be a homomorphism
on $\U(X)$.
\end{proof}

\vspace{4 mm}

Note that if $\theta_e$ is 1-1, then it follows
that $\theta$ is also 1-1.  Notice also that if
the action of $A$ on $X$ is faithful (i.e. if 
$aX = 0$ implies that $a = 0$), then $\theta_e$ is 
1-1.  This follows from Corollary \ref{es3}.
If the action of $A$ on $X$ is `completely 1-faithful'
(that is, the norm of $a \in M_n(A)$ is achieved as
the supremum of the norms of the action of $a$ on $X$
in a sensible way),
then $\theta_e$ and $\theta$ in Corollary
\ref{min} are complete isometries.

\begin{proposition}
\label{envo}  For an operator bimodule $X$, with the 
notations above, we have $I(\U_e(X)) = I(\Sy(X))$ .
If $X$ is an operator $A-B$-bimodule over C$^*-$algebras 
$A$ and $B$, which is faithful as a left and as a right 
module, and if $\U(X)$ is the triangular
operator algebra associated with a
CES representation of $X$, then $I(\U(X)) = I(\Sy(X))$ .
\end{proposition}

\begin{proof}  Clearly $\U_e(X) \subset I(\Sy(X))$.
Any minimal $\U_e(X)-$projection \cite{Ham} is the 
identity on $\Sy(X)$, and is therefore the identity map.
The rest is clear.
\end{proof}

\begin{definition} Let $Y$ be an operator $A-B$-bimodule. 
We shall call $Y$ an
{\em A-B-injective bimodule} provided for every 
pair of operator $A-B$-bimodules
V and W, with V a submodule of W, each completely 
contractive $A-B$-bimodule
map $T : V \rightarrow Y$ extends to a completely 
contractive $A-B$-bimodule map from W to Y.
\end{definition}

We should remark that the above definition corresponds to what was 
called in \cite{FP} a {\em tight} $A-B$-injective bimodule.

The following result extends Wittstock's theorem
\cite{Wi} that an injective C$^*$-algebra  is an
injective operator module over a unital C*-subalgebra.  In the 
following we consider unital C$^*-$algebras, but it is not
difficult to remove the `unital' hypothesis.

\begin{corollary} 
\label{en} Let $A$ and $B$ be unital C$^*-$algebras.
\begin{itemize}
\item [(i)]  If $Y$ is an
operator space which is also an
operator $A-B$-bimodule, then $Y$ is injective as an
operator space if and only if $Y$ is an $A-B$-injective bimodule.
\item [(ii)]  If $Y$ is an operator $A-B$-bimodule,
then the operator
space injective envelope $I(Y)$, is the operator $A-B$-bimodule
injective envelope of $Y$.  That is, $I(Y)$ is an $A-B$-injective
bimodule which is rigid and 
essential, as an operator $A-B$-bimodule containing $Y$.
Rigidity here, for example, means: any completely contractive
$A-B$-bimodule map $I(Y) \rightarrow I(Y)$, 
which is the identity on $Y$, is the identity 
on $I(Y)$.  
\end{itemize}    
\end{corollary}

\begin{proof}  One direction of (i) is obvious.  Namely,
suppose that  $Y$ is $A-B$-injective.   By CES, $Y$ may be 
realized as an $A-B$-submodule of some $B(K,H)$, where 
$H$ is a Hilbert $A$-module and $K$ is a Hilbert $B$-module. 
 Indeed $B(K,H)$ is an
operator $A-B$-bimodule.  By the $A-B$-injectivity of 
$Y$, there is a completely contractive 
projection from $B(K,H)$ onto $Y$.  Since $B(K,H)$ is
 injective as an operator space, so is $Y$.

The other direction of (i) is harder.   In a 
previous version of this 
paper we had a proof which used almost all the results 
established until now.  Instead we shall only use a few 
results from part 1, and \ref{nonchar} above.  We will also 
use the fact, which is a simple consequence of Wittstock's
original result, or Suen's modification of this result
\cite{Su}, that if $H$ is a Hilbert $A$-module, then for
any Hilbert space $K$, we have that $B(K,H)$ is $A$-injective.
Indeed their result gives
 that $B(H)$ is $A$-injective if $A \subset B(H)$ isometrically.
However in the contrary case, one may use the following kind 
of trick:  Pick a Hilbert 
space $H'$ in which $A$ is faithfully represented, then
one obtains a faithful representation of $A$ on $H \oplus H'$.
Then one may apply the Wittstock or Suen result to 
conclude that $B(H \oplus H')$ is $A$-injective, from which 
it is easy to see by compression that $B(H)$  
is $A$-injective.  See \cite{BOMD} Theorem 4.1 for
another proof of this simple consequence. 
 
Suppose that $Y$ is injective.  
Represent the C$^*-$algebra $I(\Sy(Y))$ faithfully and
nondegenerately on a Hilbert space; then the two diagonal 
projections determine a splitting of the Hilbert space
as $H \oplus K$, say.   So $I_{11}$ is a unital *-subalgebra
of B(H), and so on.  Now by injectivity there is a completely
positive projection $\phi$ 
from $B(H \oplus K)$ onto $I(\Sy(Y))$.  As in the proof of
\ref{op}, the Choi Lemma implies that this 
projection is an I(\Sy(Y))-module map, and that 
$\phi$ decomposes as a $2 \times 2$ matrix of maps.  Let
$\psi$ be the `1-2 corner' of $\phi$.  Thus
$\psi : B(K,H) \rightarrow Y$ is a completely contractive 
 projection onto $Y$, 
and its easy to see, as in \ref{op}, 
that $\psi$ is a left $I_{11}$-module map.  However if $Y$ is an
operator $A-B-$bimodule, then by \ref{nonchar} above
 there is a unital
*-homomorphism  $\theta : A  \rightarrow I_{11} 
\subset B(H)$ 
implementing the left module action.  Hence 
$H$ is a Hilbert $A$-module, via $\theta$.  Since 
$\theta$ maps into $I_{11}$, we see that the projection
$\psi$ is a left $A$-module map onto $Y$.  
Similarly, $\psi$ is a right $B$-module map onto $Y$.
Since 
$B(K,H)$ is $A-B$-injective, so is $Y$.  

(ii)  is obvious, given (i) and the fact, observed earlier, that
the $A-B$-action on $Y$ extends to make $I(Y)$ an
operator $A-B$-bimodule.
\end{proof}

\vspace{3 mm}

The proof in fact shows that any injective
operator space $X$ is
an $IM^*_l(X)-IM_r^*(X)$-injective operator bimodule.

At this point we may give another 
proof of the self-duality 
of $I(X)$:
By \ref{en} (i), $I(X)$
is injective in the category of  right C$^*-$modules over $I_{22}$.
Now appeal to \cite{Lin2} Prop. 3.10.

\vspace{3 mm}

The following corollary of \ref{co3}
generalizes a standard fact for Hilbert
spaces:

\begin{corollary}
\label{sd}  Let $Y$ be a right 
C$^*-$module over a C$^*-$algebra $A$,
which (with respect to
its canonical operator space structure) is an injective
operator space.  Then
\begin{itemize}
\item [(i)]   $Y$ is a self-dual C$^*-$module over $A$.
\item [(ii)]  $B_A(Y) = \be_A(Y)$.
\item [(iii)]  Every bounded module map $Y \rightarrow
Z$ is adjointable, for any other C$^*-$module $Z$.
\item [(iv)]  $\be_A(Y)$ is an injective C$^*-$algebra.
\end{itemize}
\end{corollary}

\begin{proof}  We may suppose w.l.o.g. that $Y$ is a full
C$^*-$module over $A$.
Thus we may regard $Y$ as a
$\K(Y)-A$-imprimitivity bimodule.
Let $Z = I(Y) = Y$ equipped with its
$\C(Y)-\D(Y)$-imprimitivity bimodule structure.  Here
$\C(Y) \subset I_{11}$ as usual.   From \cite{Ham3,Bsh},
we know that $Z \cong Y$, as imprimitivity bimodules.
Thus $\be_A(Y) \cong I_{11}$, giving (iv).
Similarly we get (i).  It is well known, and fairly
obvious,  that
(i) implies (ii) and (iii).   \end{proof}

\vspace{4 mm}

The following result must be well known.  Since we cannot
give a precise reference, we prove it:

\begin{proposition}
\label{inof}  If $A$ is a C$^*-$algebra without identity,
then the operator space injective envelope $I(A)$ is a
unital C$^*-$algebra.  Indeed $I(A) = I(A^1)$, where 
$A^1$ is the unitization of $A$.
\end{proposition}

\begin{proof}  
Represent $A$ nondegenerately and faithfully
on a Hilbert space $H$.  If $\{ e_i \}$ is an increasing
contractive approximate identity
 for $A$, then $e_i \rightarrow I_H$ in the SOT.
We also have $I_H \in A^1 \subset B(H)$.
Let $\Phi : B(H) \rightarrow I(A) \subset B(H)$ be a 
minimal $A$-projection.  So $\Phi$ is a completely contractive
idempotent map whose range contains $A$.
If $\Phi(I) = I$ we would be done, since in that case 
$\Phi$ is completely positive,
and then we can deduce that $I(A)$ is a
unital C$^*-$algebra as in the introduction (i.e. by the
Choi-Effros result quoted there).
 
In order to see that $\Phi(I) = I$,
we choose $\zeta \in H$ with $\Vert \zeta \Vert = 1$.  Let
$\phi(T) = <\Phi(T) \zeta , \zeta>$, for $T \in B(H)$.
Then $\Vert \phi \Vert \leq 1$ and $\phi(e_i) \rightarrow
1$.  It is no doubt well known and easy to see that
this implies that $\phi(I) = 1$.
One way to do this is to
write $\phi(T) = < \pi(T) \eta , \xi >$,  where $\pi$
is a unital *-representation of $B(H)$ on a Hilbert space 
$K$, and where $\eta,\xi \in BALL(K)$. 
Set $K' = [\pi(A) K]^{\bar{}}$, and let $P$ be the projection
onto $K'$.  Then it is easy to see that the net
$\pi(e_i)$  has $P$ as a WOT limit point.   For if $T$ is a
WOT limit point of $\pi(e_i)$, then for $x,y \in H$ we have   
$$<T x , y > \; = \; \lim_i < \pi(e_i) x , y > \;
= \; \lim_i < \pi(e_i) x , P y > \; = \; <  x , P y > 
\; = \; < P x , y > \; \; .$$
Thus $< P \eta , \xi > \; = \; 
\lim \phi(e_i) \;  = \;  1$.  By the converse to Cauchy-Schwarz,
$\xi = P \eta = P^2 \eta = P \xi$.  Thus 
$$\phi(I) \; = \; < \eta , \xi > \;  = \; 
< \eta , P \xi > \; = \; < P  \eta , \xi > \; = \;
1 \; . $$   
\end{proof}   

\vspace{4 mm}

In \cite{FP}, it is shown that for a C$^*-$algebra 
$A$, there is a canonical inclusion $LM(A) \subset I(A)$.
Indeed, for any essential ideal $K$ of $A$, we have 
$LM(K) \subset I(A)$.    We can offer another proof of these 
results using our methods.  Firstly, we just saw that
 if  
$X = A$ considered as an operator space, 
then $I(X)$ is a unital injective 
C$^*-$algebra.  From this it follows from abstract principles 
that $I(\Sy(X)) = I(\Sy(I(X))) = M_2(I(A))$ (see for example
\cite{Bsh} 4.14 (i)).  Thus in this case
$I_{11} = I(A)$.   Hence    
$$LM(A) = M_l(A) = IM_l(A) \subset I(A) \; . $$    

If $K$ is an essential ideal in $A$, then as we just proved,
$M(K) \subset I(K) $.   Thus we have 
$$ K \subset A \subset M(K) \subset I(K) \; .  $$
Any minimal $A$-projection on $I(K)$, is a 
$K$-projection,
and is therefore equal to the identity map, by rigidity of
$I(K)$.  Thus 
$I(K) = I(A)$.  By the first part,
$LM(K) \subset I(K) = I(A)$.
This gives the result we need.   

Notice from the above, that since $A$ is an essential ideal in 
$M(A)$, we have $I(A) = I(M(A)) = I(A^1)$.

The same is true even if $A$ is a nonselfadjoint operator algebra 
with contractive approximate identity.  In \S 4 of \cite{Bsh}
we showed that for such $A$, we have $M_l(A)$ is
the usual left multiplier algebra $LM(A)$ of $A$; and also 
that $I(A) = I(C^*_e(A))$, where 
$C^*_e(A)$ is the C$^*-$subalgebra generated by 
$A$ within the C$^*-$envelope 
$C^*_e(A^1)$ of the unitization of $A$.  It is quite clear that
$I(A^1) = I(C^*_e(A^1))$.  Hence by the above, 
$I(A) = I(C^*_e(A)) = I(C^*_e(A^1)) = I(A^1)$, which is an injective
C$^*-$algebra.  By the `abstract
principles' used a few paragraphs back, we obtain that 
$I(\Sy(A)) = M_2(I(A))$, so that $LM(A) =
IM_l(A) = \{ T \in 
I(A) : T A \subset A \}$.   Similarly
$M(A) = \{ T \in I(A) : T A \subset A , A T \subset A \}$.
Clearly any $M(A)$-projection or $LM(A)$-projection 
on $I(A)$ is an $A$-projection,
and is consequently the identity.  Thus $I(M(A)) = I(LM(A)) = I(A)$.

\vspace{4 mm}

If $X$ is a C$^*-$module we see that from results 
above that
$I(\be(X)) = I(M(\K(X))) = I(\K(X))$, where $\K(X)$ is the 
so-called `imprimitivity C$^*-$algebra' of so-called `compact' 
maps on $X$.  Does this equal $I_{11}$ in 
this case?  Certainly $\K(X) \subset I_{11}$ as a 
*-subalgebra, so that $I(\K(X)) \subset I_{11}$.

\end{document}